\input amssym.tex
\parindent=0cm 
\hoffset=1truecm
\hsize=13truecm
\vsize=23.5truecm
\font\eightrm=cmr8

\def\litem{\par\noindent\hangindent=\parindent\ltextindent}
\def\ltextindent#1{\hbox to \hangindent{#1\hss}\ignorespaces}
\long\def\ignore#1\recognize{}
\long\def\<#1\>{\vbox{\leftskip=1truecm{\eightrm #1}}} 
\def\big{\bigskip}
\def\med{\medskip}
\def\meds{\medskip}
\def\cl{\centerline}
\def\ol{\overline}
\def\sm{\setminus}
\def\ld{\ldots}
\def\cd{\cdot}
\def\to{,\ldots,}

\def\{{\lbrace}
\def\}{\rbrace}

\def\map{\rightarrow}
\def\inv{^{-1}}
\def\N{{\Bbb N}}

\def\Q{{\Bbb Q}}
\def\coeff{{\rm coeff}}

\def\label{{\rm lab}}

\def\weak{{\curlyvee}}
\def\prior{{\circ}}
\def\short{{\circ}}
\def\ord{{\rm ord}}
\def\min{{\rm min}}

\def\inter{{\sim}}
\def\J{{\cal J}}
\def\I{{\cal I}}
\def\M{{\cal M}}
\def\K{{\bf K}}

\def\K{{\cal K}}

\def\a{\alpha}
\def\mm{m}
\def\abs#1{\vert#1\vert}
\def\nc{normal\ crossings\ }
\def\usc{upper\ semicontinuous\ }

\def\ttop{{\rm top}}

\def\strat{{\rm strat}}


\cl {\bf Strong resolution of singularities in characteristic zero}\medskip\bigskip

\cl {Santiago Encinas and Herwig Hauser}\bigskip\bigskip


Hironaka's spectacular proof of resolution of singularities is built on a multiple
and intricate induction argument. It is so involved that only few people could really
understand it. The constructive proofs given later by Villamayor, Bierstone-Milman
and Encinas-Villamayor presented important steps towards a better understanding of the
reasoning. They describe an algorithmic procedure for resolution, using a local
invariant to show that the situation improves under blowup. The centers of blowup are
given as the locus where the invariant takes its maximal value. \med

Despite this simple outset, the proofs are not easy, mostly because the
definition of the invariant requires to know the whole sequence of blow\-ups having
occurred so far in the resolution process. To make the invariant intrinsic and to patch
the local constructions various equivalence relations have to be introduced.\med

When trying to understand these proofs it became clear that to define the
invariant it is actually not necessary to refer to the entire sequence of earlier
blowups. It is sufficient to have information on two sets of exceptional components at
each stage of the process. Including this information to the resolution datum (called
{\it mobile} in this paper) the invariant can be defined directly without going back to
the very beginning of the resolution process.  \med

{\it ``If you wish to cross the Sahara, better take a map with you where you mark 
daily the trajectory you have made so far rather than to return every morning to
the starting point in order to know in which direction to continue.''} \med

Already Abhyankar and Hironaka took a map with them, and the idea is also used in the
other papers, though in a different way. The extra luggage we shall carry with us
specifies exactly the combinatorial book-keeping information which allows to
construct the local invariant. The resulting proof should be regarded as a conceptualized
version of the existing proofs. In substance it is the {\it same} proof.\med

Although the invariant can now be defined promptly potential readers may still
hesitate to try to {\it understand} how to prove resolution of singularities. For
many,  it is a black box better not to be touched. But Hironaka's proof is a
phantastic proof, and very beautiful. Our objective is to transmit this fascination --
and to open the box. So let us start.\big\med\goodbreak

Keywords: Resolution of singularities.

AMS Classification: 14E15, 14B05, 32S05, 32S10, 32S45
\big\med\goodbreak


{\bf The result}\med

Let $X$ be a reduced singular scheme. A {\it strong resolution of $X$} is,
for every closed embedding of $X$ into a regular ambient scheme $W$, a proper
birational morphism $\varepsilon$ from a regular scheme $W'$ onto $W$ subject to the
following conditions.\medskip

{\leftskip=.6truecm

{\it Explicitness}. $\varepsilon$ is a composition of blowups of $W$ in regular closed
centers $Z$ transversal to the exceptional loci.\par

{\it Embeddedness}. The strict transform $X'$ of $X$ is regular and has normal crossings 
with the exceptional locus in $W'$.\par

{\it Excision}. The morphism $X'\map X$ does not depend on the embedding of $X$ in
$W$.\par
 
{\it Equivariance}.  $\varepsilon$ commutes with smooth morphisms $W^-\map W$,
embeddings $W\map W^+$, and separable field extensions. In particular, $\varepsilon$ is
an isomorphism outside the singular locus of $X$ and commutes with group actions. \med

}

\<The resolution of $X$ is global with centers  equal to the top locus of an \usc invariant
$i_a(X)$ on $W$ given by the local rings of $X$ ({\it effectiveness}).  The resolution
commutes with open immersions, local and global diffeomorphisms and taking cartesian products
with regular schemes. The smooth morphisms of {\it equivariance} need not be defined over the
ground field. Passage to the completions implies resolution of formal schemes. The resolution
process can be implemented.  \>\meds


The existence of resolutions satisfying the first two properties was established for
excellent schemes of finite type over a field of characteristic zero by Hironaka [H1].
 For such schemes, the construction of a strong resolution in the above sense through a
local invariant defining the centers is due to Villamayor [V1, V2], Bierstone-Milman [BM1,
BM2, BM3] and Encinas-Villamayor [EV1, EV2]. The algorithm was implemented by
Bodn\'ar-Schicho [BS].  \medskip

Weak resolution theorems in characteristic zero have been established with 
different methods by Abramovich-de Jong, Abramovich-Wang and Bogomolov-Pantev [AJ, AW, BP].
The main results on resolution in positive characteristic are due to Abhyankar, 
Giraud, Lipman, Cossart, de Jong and Kuhlmann. We refer to [Ha1] for details. \medskip 

The proof of the present paper for the existence and construction of strong resolutions
in characteristic zero relies on ideas and techniques from Hironaka, Abhyankar and
Giraud. The invariant and the centers are almost identical to the ones used by
Vil\-la\-mayor, Bierstone-Milman and Encinas-Villamayor. There are, however, several
improvements with respect to the existing literature. \medskip


The resolution of schemes will be deduced from the resolution of {\it singular mobiles}.
Mobiles are intrinsic global data which encode the singular and combinatorial structure of 
the resolution problem and record its transversality with the exceptional divisor. Locally,
each mobile is exploited through the choice of a {\it punctual setup} associated to it.
This is a  string of ideals in decreasing dimensions which determines the resolution
invariant. The latter is shown not to depend on the chosen setup. It thus gives a local
measure of the resolution complexity of the mobile. Its {\it top locus} defines a global
center of blowup, which in turn determines the {\it transform} of the mobile. The 
invariant of the transformed mobile is shown to have decreased and thus induction
applies to give the resolution of mobiles and hence of schemes. \med

The main issues of the proof are the following.\med 


{\leftskip=.6truecm

Coverings of global and patchings of local objects are avoided by considering {\it
stratified} ideals and punctual setups of mobiles defined only {\it locally}. This and the 
use of {\it transversal handicaps} eliminates also the use of equivalence relations.\med

All constructions and arguments but one are {\it characteristic free}, the exception being
the existence of osculating hypersurfaces.\med

The Hilbert-Samuel function and normal flatness are avoided by working with
the order of ideals and weak transforms. The centers of blowup lie in the weak transforms
of the scheme though possibly not in the strict transform. This is not a serious drawback,
since all centers map to the singular locus, but simplifies things considerably. For
equidimensional schemes, the Hilbert-Samuel function had already been eliminated by
Encinas-Villamayor in [EV2] using a somewhat different argument; for the general case, 
see also [BV]. \med

The relevant information on the history of the resolution process is encoded in the
{\it combinatorial} and {\it transversal handicap} of the mobile and its transforms. To define
the invariant it is no longer necessary to consider the whole sequence of earlier blowups.\med

}


Standard results on the order of ideals and on hypersurfaces of maximal contact as well
as some straightforward verifications are omitted.  Paragraphs in small characters
provide background information and/or proofs of well known or technical results
appearing in the main body of the text.  Superscripts  refer to the appendix. For
notational convenience there appear rational powers of ideals. These could be avoided by
raising the ideals to suitable powers,  as will be indicated in parentheses.  As taking
the order commutes with powers, the exponents are treated as integers.\med


The various constructions of the present paper are often justified only {\it a
posteriori} through their role in the proof. This makes it hard to get a feeling for
them at the moment when they are introduced. The expository paper [Ha2] provides these
justifications {\it ab initio}. It shows how the constructions arise naturally when
trying to prove resolution of singularities from scratch. \med


We are indebted to Hironobu Maeda for very valuable references, and to G\'abor Bodn\'ar for
many probing questions.  Substantial improvements regarding the organization of the
paper were suggested by a highly competent anonymous referee.  \medskip

\bigskip\goodbreak


{\bf Idea of proof}\med 

Let $J$ be the ideal  defining $X$ in $W$. We want  to transform $J$ by a sequence
of blowups into a simple form, i.e., so that the pull-back of $J$  becomes a monomial principal
ideal.  The resolution of $X$ will be then deduced from this monomialization.  Let us
place ourselves at a certain stage of the resolution process. We will have to decide on
the center of the next blowup. The ideal $J$ will stem from earlier blowups, so that
exceptional components can be factored from $J$ to a certain power. This factor will be
noted down in what we call the combinatorial handicap. It is a (non-reduced) normal
crossings divisor $D$ in $W$ supported by the current exceptional locus $F$ so that $J$
factors into $J=M\cd I$ with $M$ the ideal defining $D$ in $W$, and some ideal $I$ of
$W$ which is still unresolved.  Our objective will be to lower the order $o$ of $I$ 
at the points of  $W$ by further blowups, until $I$ becomes $1$ and $J=M$ is the
required monomial. A separate argument will show that the monomialization of ideals
implies the resolution of singular schemes.\med

Fix the above situation. The center $Z$ of the next blowup $\pi:W'\map W$ should be a
closed and globally defined regular subscheme of $W$, which is transversal to the
exceptional locus. In addition, we wish to have $Z$ inside the top locus of $I$, i.e., in
the set of points where the order of $I$ in $W$ is maximal. In particular,
$o=\ord_aI=\ord_ZI$ shall hold for all $a\in Z$. Here, $\ord_ZI$ denotes the maximal power of
the ideal of $Z$ in $W$ which contains $I$. This will ensure that the order $o'$ of the
transform of $I$ under the blowup of $W$ with center $Z$ will not increase. Once $Z$ satisfies
these two conditions, we will have $o'\leq o$ for all points of the new exceptional component
$Y' =\pi\inv(Z)$ in $W'$, and the total transform of $M$ will be an ideal $M^*$ defining again
a normal crossings divisor in $W'$. By construction, the total transform $J^*$ of $J$ will
factor into $J^*=M^*\cd I(Y')^o\cd I^\weak$, where $I^\weak=I(Y')^{-o}\cd I^*$ denotes the
weak transform  of $I$. Setting $J'=J^*$, $M'= M^*\cd I(Y')^o$ and $I'=I^\weak$ we get again a
product \med

\cl{$ J'=M'\cd I'$}\med

with prescribed exceptional factor $M'=I_{W'}(D')$ given by the transformed
combinatorial handicap $D'=D^* + o\cd Y'$, and new exceptional locus
$F'$. Thus our (preliminary) resolution datum, made precise later through the concept of
mobiles, consists at each stage of a product of ideals $J=M\cd I$ and two normal
crossings divisors $D$ and $F$ in $W$.\med

We are left to determine a suitable center $Z$, and to show that at
the points where equality $o'=o$ holds the situation has improved.  Both tasks will
be accomplished simultaneously by associating to $J$, $D$ and $F$ a local upper
semicontinuous invariant $i_a(J)$. Its top locus will be the required center $Z$,
and $i_a(J)$ will drop after blowup. \med

The crucial advantage in characteristic zero is that there exists locally at each
point a regular hypersurface $V$ of $W$ whose successive transforms under any blowup
with center inside $V$ contain all equiconstant points, i.e., the points where the orders
$o'$, $o''$, ... of the transforms of $I$ remain constant (hypersurface of maximal
contact). Choose such a $V$ locally at a point $a$ of $W$. Let $Z$ be any center of
blowup inside $V$. As the transform $V'$ of $V$ contains the points $a'$ above $a$ where
$o'=o$ it suffices to compare $J$ and $J'$ at points of $V$ and $V'$, i.e., inside
hypersurfaces.\med

The idea then is to associate to $J=M\cd I$ and $J'=M'\cd I'$ ideals $J_-$ and
$(J')_-$ in $V$ and $V'$ which reveal the expected improvement. Once we have constructed
the appropriate ideal $J_-$ we can apply induction on the dimension to find the center and
the invariant, since $J_-$ lives in an ambient space of smaller dimension. In this way we
may assume that we have already constructed a local upper semicontinuous invariant
$i_a(J_-)$ whose top locus prescribes a regular center $Z_-$ in $V$    such that
blowing up $V$ in $Z_-$ makes $i_a(J_-)$ drop (except if $J_-$ is already resolved). If
$i_a(J_-)$ does not depend on any choices (in particular, not on the local choice of
$V$), the center will automatically be defined globally. We will give the definition of
$J_-$ in a moment.\med

There arise two problems. The center $Z_-$ associated to $J_-$ may not be transversal 
to the exceptional locus $F$,  and the transform of $J_-$ under the blowup $V'\map V$
with center $Z_-$ may not coincide with the ideal $(J')_-$ associated to $J'$ in $V'$.\med

If $Z_-$ is not transversal to $F$, we have to solve this subproblem first. Auxiliary
blowups with smaller centers will make $Z_-$ transversal to $F$, so that it can be really
taken as center. Actually, $J_-$ will be built up so that this subproblem is solved in
parallel: we specify the components $E$ of $F$ to which $Z_-$ may not be transversal, 
noted down in the transversal handicap, and then resolve the ideal $Q=I_V(E\cap V)$ in $V$
by auxiliary blowups. Once its weak transform has become $1$, $V$ and $E$ will be
separated from each other, and transversality holds since $Z_-\subset V$. This separation
cannot and need not be realized for the whole exceptional locus $F$: the components of
$F\sm E$ will a priori be transversal to $Z_-$ and therefore do not affect the
transversality problem. Of course, the critical components $E$ inside $F$ have to be
determined explicitly.\med

The second problem is handled by taking for $(J_-)'$ an intermediate transform
bet\-ween total and weak transform (the controlled transform; it is given by a number
$c$, the control). For this, the required commutativity $(J')_-=(J_-)'$ is a 
check in local coordinates. It is here that we need to work with factorized ideals 
$J=M\cd I$, because $J$ and $I$ will transform differently. The controlled transform
of $J$ ensures commutativity while its order may increase, whereas the weak transform of $I$
would not yield commutativity while its order decreases or remains constant.\med

We see that everything concentrates on defining the correct ideal $J_-$. This
will be achieved through the coefficient ideal of $I$ in $V$. It is
obtained from $I$ by expanding its elements with respect to a local coordinate
defining $V$ in $W$ and taking the ideal generated by (equilibrated powers of) its
coefficients. To include the transversality problem, one has to take the
coefficient ideal not of $I$ but of a product $I\cd Q$ where $Q$ defines the
possibly non-transversal components from $E$. Then $Z_-$ will be
contained in the top locus of $Q$, hence in all components of $E$. Therefore it will
automatically be transversal to $F$. \med

It remains to define the invariant $i_a(J)$. It is given as the vector  \med

\cl{ $i_a(J) =(\ord_aI, \ord_aQ, m_a, i_a(J_-))$,}\med

where the component $m_a$ is of combinatorial nature. It only becomes relevant when $J_-$
is already resolved and its invariant $i_a(J_-)$ cannot improve.   The invariant is
considered with respect to the lexicographic order. It depends on $J$, $D$ and $E$.\med

If $J_-$ is not resolved, we may assume by induction on the dimension that
$i_a(J_-)$ will improve when blowing up its top locus  $Z_-$ in $V$ (the case of 
dimension of $V$ equal $1$ being trivial). Eventually, $(\ord_aI, \ord_aQ)$ must drop. When
the second component $\ord_aQ$ drops, the transversality subproblem improves. After finitely
many steps it is solved, $\ord_aQ=0$. Larger centers become permissible. As $i_a(J_-)$
continues to improve the first component $\ord_aI$ must also drop sometime. Now induction
applies to prove that finitely many blowups yield $\ord_aI=0$, which signifies that $J=M$
is a monomial.\med

The above considerations show that we need in each dimension informations on the
divisor $D$ formed by those exceptional components which can be factored from $J$ and on
the divisor $E$ of exceptional components which may pose a transversality
problem. Therefore the combinatorial and the transversal handicap appearing in mobiles
will consist of strings $D_n\to D_1$ and $E_n\to E_1$ of (stratified) normal crossings
divisors in $W$. They are global  objects which do not depend on any local choices,  and
obey prescribed rules of transformations under blowup. Thus we can define the transform of a
mobile under blowup, and its resolution. Mobiles are the maps we take with us on our trip
through the Sahara: Every day  we write down how the combinatorial and transversal handicap
have transformed under the last blowup. The transformation rule for the handicaps will depend
on the point of the blowup we are considering (i.e., on the value of the local invariant at
this point). As the invariant is upper semicontinuous and hence induces a stratification by
locally closed sets, we will get stratified divisors.\med

The descent in dimension via local ideals $J_-$  depends on the choice of the
hypersurface $V$  and yields strings of ideals $J_n\to J_1$ (with $J_n=J$, $J_{n-1}=J_-$) in
local flags $W_n\supset\cdots \supset W_1$ of regular subschemes, called the setup of the
mobile. Each ideal factors into $J_i=M_i\cd I_i$ according to $M_i=I_{W_i}(D_i \cap W_i)$, and
$J_i$ will be the coefficient ideal in $W_i$ of some ideal $K_{i+1}$ associated to
$J_{i+1}$ (the ideal $K_i$ is essentially the product $I_i\cd Q_i$ with
$Q_i=I_{W_i}(E_i\cap W_i)$ given by the transversal handicap). The ideals $I_i$ and $Q_i$
depend on choices, but their orders, which form the components of the invariant
$i_a(J)$, do not. Thus $i_a(J)$ is intrinsic, and upper semicontinuous because
orders of ideals are. \med

It then has to be shown that the top locus of $i_a(J)$ is in fact regular and
transversal to the exceptional locus.  This allows to choose it as the center of the
next blowup. It remains to prove that the transform of the mobile in $W'$ admits at
each point a setup which is the transform of the setup of the initial mobile
(commutativity), thus its invariant can be computed from the invariant below and must have
decreased. As the invariant $i_a(J)$ can only drop finitely many times, it must eventually
achieve its minimal value $0$. In this case, $\ord_aI=0$ and $J=M$ as required.\med 

The above proof is based on a cartesian scheme of induction: the descending
horizontal induction on the dimension is combined with the vertical induction on the
resolution invariant. In total, fourteen inductions are required. Once the relevant
objects like mobiles, setups and their transforms  are defined properly, the inductions
follow always the same pattern.

\big\med\goodbreak


\med \cl{\bf CONSTRUCTIONS}\bigskip

{\bf Concepts}\medskip

Throughout, we fix a regular ambient scheme $W$ and a regular locally closed $n$-dimensional
subscheme $V$ of $W$. By a {\it divisor} in $W$ we shall mean an effective Weil divisor $D$. A
closed subscheme $D$ of $W$ has {\it normal crossings} if it can be defined  locally by a
monomial ideal. The subscheme $V$ meets $D$ {\it transversally} if the product of the defining
ideals of $V$ and $D$ defines a \nc scheme. \medskip

A {\it local flag} in $V$ at $a$ is a decreasing sequence $W_n\supset\ld\supset W_1$ of closed
$i$-dimensional regular subschemes $W_i$ of a neighborhood $U$ of $a$ in $V$. An ideal $K$ in
$V$ is {\it bold regular} if it is a power of a regular principal ideal in $V$. A {\it
stratified ideal} in $V$ is a collection of coherent ideal sheaves each of them defined on a
stratum of a stratification of $V$ by locally closed subschemes.  A {\it stratified
divisor} is defined by a stratified principal ideal. All ideals and divisors will be stratified
without notice, except if said to be coherent. For \nc divisors $D$, the underlying
stratification $\strat(D)$ need not have \nc strata.  \medskip\goodbreak

A map $Q_b\map (Q_b)^\sharp$ associating to stalks of ideals $Q_b$ in an open subscheme
$U$ of $W$ stalks of ideals $(Q_b)^\sharp$ in $V$ is {\it tuned} along  the stratum $S$ of a
stratification of $V$ through a point $b$ of $V$  if $Q_b$ and $(Q_b)^\sharp$ admit
locally at any point $b$ of $V$ coherent representatives $\ol {Q_b}$ on $U$ and 
$\ol{(Q_b)^\sharp}$ on $V$ so that the stalks $((\ol {Q_b})_a)^\sharp$ and
$(\ol{(Q_b)^\sharp})_a$ at $a$ coincide along $S$.  This is abbreviated by saying that the
ideals $(Q_b)^\sharp$ are tuned along the stratum $S$.  \medskip


A {\it shortcut} of a \nc divisor $M$ in $W$ is a divisor $N$ obtained from $M$ by 
deleting on each stratum  of the underlying stratification $\strat(M)$ of $M$ some
components of $M$.  The divisor $M$ is {\it labelled} if each  shortcut $N$ comes with a
different non negative integer $\label\,N$, its {\it label}. The empty shortcut has
label $0$. A shortcut $N$ of a \nc divisor $M$ is {\it tight at $a$ of order} $\geq c$
if it has order $\geq c$ at $a$ and if any proper shortcut of $N$ has order $<c$ at $a$.
It is {\it maximal tight at} $a$ if $M$ is labelled and if $(\ord_a N,\label\,N)$ is
lexicographically maximal among the tight shortcuts of $M$ of order $\geq c$ at $a$
${}^{<1>}$. \medskip


A {\it handicap} on $W$ is a sequence $D=(D_n\to D_1)$ of stratified \nc divisors
$D_i$ of $W$. The {\it truncation} of $D$ at index $i$ is ${}^iD=(D_n\to D_i)$. \medskip


A {\it singular mobile} in $W$ is a quadruple $\M=(\J,c,D,E)$ with $\J$ a coherent nowhere
zero ideal sheaf on $V$, $c$ a  non negative constant associated to $V$  and $D$ and $E$
handicaps in $W$ with $D$ labelled  and $E$ reduced${}^{<2>}$. We call $c$ the {\it
control} of $\J$, and $D$ and $E$ the  {\it combinatorial} and {\it transversal handicap} of
$\M$. The truncation ${}^i\M$ at index $i$ of $\M$ is $(\J,c,{}^iD,{}^iE)$.\medskip\goodbreak

The {\it transversality locus} of a mobile $\M$ is  $\abs E=E_n\cup\ld\cup E_1$. The 
{\it exceptional locus} of a sequence of blowups of $\M$ is the reduced inverse image
of the union of all centers and of $\abs E$.\med


\<The control $c$ is allowed to be $0$ if and only if $J=1$. Throughout, we denote by
calligraphic letters stratified ideal sheaves, and by roman letters their stalks or
sufficiently small representatives of them.\>\meds


A {\it strong resolution of a mobile} $\M=(\J,c,D,E)$ in $W$ with $\J$ a nowhere zero ideal in 
$V$ is a sequence of blowups of $W$ in regular closed centers $Z$ such that the ideal $\J'$ of
the final transform $\M'=(\J',c',D',E')$ of $\M$ as defined in the section {\it ``Transform of
mobile''} has order $<c$. We require that the centers are transversal to the exceptional loci,
and that the resolution is equivariant.\medskip


The {\it top locus} of an upper semicontinuous function $t$ on  $V$  is the reduced closed
subscheme $\ttop(t)$ of points of $V$ where $t$ attains its maximum. The {\it order} at
$a$ of an ideal $\J$ of $V$ is the largest power $o=\ord_a\J$ of the maximal ideal of 
${\cal O}_{V,a}$ containing the stalk $J$ of $\J$ at $a$. We set $\ttop(\J)=\ttop(\ord\,\J)$
and denote by $\ttop(\J,c)$ the locus of points in $V$ where the order of $\J$ is at
least $c$.  For closed subschemes of $V$, the analogous loci are defined through the
associated ideals. When working locally at a point $a$, $\ttop(t)$ also denotes the
local top locus of $t$ in a neighborhood of $a$.  \medskip 


Let $W'\map W$ be the blowup of $W$ with center $Z$  inside $V$ and exceptional component
$Y'$. The {\it total} and {\it weak transform} of an ideal $\J$ of $V$ are the inverse 
image $\J^*$ of $\J$ under the induced blowup $V'\map V$ and the ideal $\J^\weak=\J^*\cd
\I(Y'\cap V')^{-o}$ with $o=\ord_Z\J$.  The  {\it controlled transform} of $\J$ with respect
to $c\leq o$ is the ideal $\J^!=\J^*\cd \I(Y'\cap V')^{-c}$ in $V'$ ${}^{<3>}$. \medskip

\goodbreak 


The {\it companion ideal} $P$ of a product $J=M\cd I$ of ideals in $V$ at $a$ with
respect to a control $c\leq \ord_aJ$ on $V$ is the ideal $P$ in $V$ at $a$ given by\medskip

\hskip 3cm $P=I + M^{o\over c-o}$\hskip 2cm if \ $0<o=\ord_aI< c$,\medskip

\hskip 3cm $P=I $\hskip 3.3cm otherwise.\medskip

\<To avoid rational powers of ideals, one could take $P=I^{c-o} + M^o$ if \
$0<o=\ord_aI< c$, and $P=I$ otherwise.\>\meds 

Companion ideals are tuned along  $\ttop(I)\cap \strat(M)$  and satisfy
$\ttop(P)=\ttop(I)\cap\ttop(M,c-o) \subset\ttop(J,c)$. The {\it transversality ideal} $Q$ in
$V$ of a \nc divisor $E$ of $W$ is the ideal \medskip

\hskip 3cm $Q=I_V(E\cap V)$ \medskip

defining $E\cap V$ in $V$.  The  {\it composition ideal}  $K$ in $V$ of a product $J=M\cd
I$ of ideals in $V$ at $a$ with respect to a control $c$ and a \nc divisor $E$ in $W$ equals
${}^{<4>}$\medskip
 
\hskip 3cm $K=P\cd Q$ \hskip 2cm if $I\neq 1$,\medskip

\hskip 3cm $K=1$ \hskip 2.6cm if $I= 1$,\medskip

with $P$ the companion ideal of $J$ and $c$, and $Q$ the transversality ideal of
$E$ in $V$. Composition ideals are tuned along  $\ttop(I)\cap \strat(M)\cap \strat(E)$ 
and satisfy locally $\ttop(K)= \ttop(P)\cap\ttop(Q)$ if $I\neq 1$ and $Q\neq 0$.
\medskip\goodbreak


\<The weight on $M$ in the definition of $P$ is necessary to have $P^\weak=(I + M^{o/
c-o})^\weak=I^\weak + (M')^{o'/c'-o'}$ when $c'=c$ and $o'=o$.  This ensures that the
passage to companion ideals commutes with blowups.  The equality
$\ttop(P)=\ttop(I)\cap\ttop(M,c-o)$ uses that $\ord_aM^{o/c-o}$ is $\geq o$ for $a$ in
$\ttop(J)$ since there $c\leq \ord_aJ=o+\ord_aM$. Note that $\ord_aP=\ord_aI=o$ along
$\ttop(P)\subset\ttop(I)$, and thus $P=1$ if and only if $I=1$.\>\meds
 
\<In the application, $V$ will be a member $W_i$ of the local flag and $E$ will be a
member $E_i$ of the transversal handicap; the latter does not contain
$W_i$ and thus $Q_i\neq 0$. Moreover, $W_i$ will be transversal to $E_i$ and thus
$\ord_aQ_i=\ord_aE_i$, i.e., equal the number of components of $E_i$ passing through $a$,
since $E_i$ is reduced. We have $\ttop(Q_i)=W_i\cap \ttop(E_i)$ locally at $a$. Therefore
any $Z\subset\ttop(Q_i)$ is contained in the components of $E_i$ it meets.\>\meds 


The {\it tag} of an ideal $J$ in $V$ at $a$ with control $c$ and \nc divisors $D$
and $E$ in $W$ such that $J=M\cd I$ for $M=I_V(D\cap V)$ with $D$ labelled  and transversal
to $V$  is the vector  \medskip

\hskip 3cm $t_a(J)=(o,k,\mm)\in\N^4$, \medskip\goodbreak

equipped with the lexicographic order. Here, $o=\ord_aI$ and $k=\ord_aK$ with $K=P\cdot
Q$ the composition ideal of $(J,c,E,D)$. We set $\mm=(0,0)$ if $o>0$, and
$\mm=(\ord_aN,\label\,N)$ otherwise with $N$ the maximal tight shortcut of $M$ at $a$ of
order $\geq c$. \medskip
 

\<The combinatorial tag $\mm$ measures the improvement of the controlled transform of $M$
once we have $J=M$, the singular tag $(o,k)$ measures the improvement of the weak transforms
of $I$ and $K$.\> \meds


The {\it coefficient ideal} ${}^{<5>}$ of an ideal $K$ of $W$ at $a$ with respect to 
$V$ is an ideal  in $V$ which is built from the coefficients of the Taylor expansion of
the elements of $K$ with respect to the equations defining $V$. Let $x,y$  and $y$ be
regular systems of parameters of ${\cal O}_{W,a}$ and ${\cal O}_{V,a}$ so that $x=0$
defines $V$ in $W$. For $f$ in $K$ denote by $a_{f,\alpha}$ the elements of ${\cal
O}_{V,a}$ so that $f= \sum_\alpha a_{f,\alpha}\cd x^\alpha$ holds after passage to the
completion.  Then set \medskip

\hskip 3cm $\coeff_VK=\sum_{|\alpha|<c}(a_{f,\alpha}, f\in
K)^{c\over c-|\alpha|}$.\medskip
 
 Replacing the exponents by ${c!\over c-|\alpha|}$, rational powers of ideals could be
avoided.  Coefficient ideals are tuned along $\ttop(K)\cap V$. Let $V'\map V$ be the
blowup of $V$ induced by $W'\map W$ with center $Z$ contained in $\ttop(K)\cap V$ and
exceptional component $Y'\cap V'$. At points $a'$ of $Y'\cap V'$ where
$c'=\ord_{a'}K^\weak=\ord_aK=c$ one has $ (\coeff_VK)^!=\coeff_{V'}(K^\weak)$. Thus the
descent to the coefficient ideal commutes with taking the weak, respectively controlled
transform at these points.\medskip 


\<The proof of commutativity goes as follows, using that
$a_{f^{\weak},\alpha}=(a_{f,\alpha})^{*}\cdot I(Y'\cap V')^{|\alpha|-c}$:\>\meds 


\<\hskip 2cm $\coeff_{V'}(K')=\coeff_{V^\weak} (K^\weak) =$\med

\hskip 2cm $=\coeff_{V^\weak} (\sum_{\abs\a<c'} a_{f^\weak,\a}\cd x^\a,\ f^\weak\in
K^\weak)=$\med

\hskip 2cm $ =\coeff_{V^\weak} (\sum_{\abs\a<c'} a_{I(Y'\cap V')^{-c}\cd f^*,\a}\cd
x^\a,\ f^\weak\in K^\weak)=$\med

\hskip 2cm $ =\coeff_{V^\weak} (\sum_{\abs\a<c} (a_{f,\a }\cd x^{\abs{\a}})^*\cd
I(Y'\cap V')^{-c},\ f\in K)=$\med

\hskip 2cm $=\sum_{\abs\a<c} (a_{f,\a}^*,\ f\in K)^{c/(c-\abs\a)}\cd
I(Y'\cap V')^{-c}=$\med

\hskip 2cm $=(\sum_{\abs\a<c} (a_{f,\a},\ f\in K)^{c/(c-\abs\a)})^*\cd
I(Y'\cap V')^{-c}=$\med

\hskip 2cm $=I(Y'\cap V')^{-c}\cd (\coeff_VK)^*= (\coeff_VK)^!$.\>\meds


\<The coefficient ideal depends on the choice of $V$ and the regular systems of
parameters. It satisfies $\ord_a\coeff_VK\geq \ord_aK$ for $a\in V$.\>\meds


The {\it junior ideal} $J$ in $V$ of an ideal $K$ of $W$ at $a$ is the coefficient ideal
$\coeff_VK$ of $K$ in $V$ if $K$ is not bold regular or $1$, and is set equal to $1$ otherwise.
\medskip

The scheme $V$ has {\it weak maximal contact} with an ideal $P$ of $W$ at $a$ if $V$
maximizes the order of $\coeff_VP$ at $a$. It is {\it osculating} ${}^{<6>}$ for $P$ if 
there is an $f\in P$ with $\ord_af=\ord_aP$ and $\ord_a\coeff_Vf= \ord_a\coeff_VP$ such that
$a_{f,\alpha}=0$ for all $\alpha$ with $|\alpha|=\ord_aP-1$. \medskip


\<The coefficient ideal of $K\neq 0$ in a hypersurface of weak  maximal contact $V$
is zero if and only if $K$ is bold regular or equal $1$. Namely, if $K$ is bold regular then
$\ttop(K)=V$, by weak maximal contact, and the coefficient ideal is zero. Conversely, if
$K\neq 1$ and the coefficient ideal is zero, then  $a_{f,\alpha}=0$ for $f\in K$ and
$\vert\alpha\vert<c$. Thus $K\subset(x^c)$ with $x$ defining $V$ in $W$. Hence the support
of $K$ contains the hypersurface defined by $x^c$. As the order of $K$ is at most $c$ we
get $K=(x^c)$ bold regular. \>\meds


\<Junior ideals avoid getting zero ideals $J$ and $I$ having infinite order. They are tuned
along  $\ttop(K)\cap V$.  The passage to junior ideals in hypersurfaces  of weak maximal
contact commutes with taking the weak, respectively controlled transform, at points $a'$ above
$a$ where the order of $K$ has remained constant. This holds for coefficient ideals. If $K$ is
not bold regular or $1$ but $K'$ is this implies that the order of $K$ has dropped, by the
commutation of coefficient ideals with weak and controlled transforms. If $K$ is bold
regular but not $1$  commutativity does not hold  except if the center equals the support
of $K$. But by definition of $Z$, it does equal the support of $K$ in this case.\>\meds


 The following properties of osculating hypersurfaces will be used later (see below for
proofs.)  In characteristic zero, any $a$ in $W$ admits a neighborhood $U$ and a
hypersurface $V$ of $U$ which is osculating for $P$ at all points of $\ttop(P)$ in $U$.
If $P\neq 0,1$  then $V$ contains $\ttop(P)\cap U$ and satisfies
$\ttop(J_P,c_P)=\ttop(P)$ for $c_P=\ord_aP$ and $J_P$ the coefficient ideal of $P$ in
$V$. Any $V$ which is osculating for $P$ has weak maximal contact with $P$.  If  $V$
has weak maximal contact with $P\neq 1$ then it has also weak maximal contact with any
product $K=P\cd Q$. \meds

If $V$ has weak maximal contact with $P$ its weak transform $V^\weak$ under blowup of $W$ in 
a  regular center $Z\subset \ttop(P)$ contains all points $a'$ of $W'$ where the order of
$P^\weak$ has remained constant. If $V$ is osculating for $P$, then $V^\weak$ is osculating
for $P^\weak$ at these points. Hence, if $V$ has weak maximal contact with $K=P\cd Q$ being
osculating for $P$ and if $Z\subset\ttop(K)$ then $V^\weak$ has weak maximal contact with
$K^\weak=P^\weak\cd Q^\weak$ at all points $a'$ of $W'$ where $\ord_{a'}P^\weak=\ord_aP$,
regardless of the order of $K^\weak$ at $a'$. \medskip 


\< We indicate how to prove the more delicate of the above assertions. For the
existence of osculating hypersurfaces, take an element $f\in P$ of minimal order
$\ord_af=\ord_aP$ at $a$. In characteristic $0$, a suitable partial derivative of $f$ of 
order $\ord_af-1$ has order $1$ at $a$ and can be taken as defining equation for $V$.
From this, $\ttop(J_P,c_P)=\ttop(P)\subset V$ is immediate. Equivalently, $V$ can be
obtained by a Tschirnhaus coordinate transformation${}^{<6>}$.  \>\meds


\<That osculating with $P$ implies weak maximal contact with $P$ is a direct check on 
Newton polyhedra.  The properties with respect to blowup are classical and proven by
computing the transforms of the ideals in local coordinates for which the blowup
reduces to a monomial substitution of variables.  We show that weak maximal contact
with $P\neq 0, 1$ implies weak maximal contact with $K=P\cd Q$. Let $c$, $c_P$ and $c_Q$
be the orders of $K$, $P$ and $Q$. We may assume that $Q\neq 0, 1$, so that all orders
are positive and finite. Let $J$, $J_P$ and $J_Q$ be the associated coefficient ideals,
and $e$, $e_P$ and $e_Q$ their respective orders in $\N\cup\{\infty\}$.  We claim that
\medskip

\> \<       \hskip 3cm $e/c=\min\{e_P/c_P,e_Q/c_Q\}$. \medskip

The equality also holds when $P$ is bold regular, in which case $J_P=0$ and $e_P=\infty$.\>\meds

\<To prove the claim, we fix local coordinates $(x,y_{n-1}\to y_1)$ at $a$ in $W$ with
$V$ defined in $W$ by $x=0$. We treat first the  case of single elements $f\in P$ and 
$g\in Q$. Their Newton polyhedra satisfy $N(f\cd g)=N(f)+N(g)$, as is shown by a
computation of the vertices of $N(f\cd g)$. Interpret then the order $e_{f\cd g}$  of
the ideal generated by the equilibrated powers of the coefficients of $f\cd g$ in $J$
 as the order of the projected Newton polyhedron $\pi(N(f\cd g))$ under the projection
$\pi$ from the point $(c,0\to 0)$ of $\N^n=\N\times \N^{n-1}$ to $\Q^{n-1}$ defined for
points $(j,\gamma)$ with $j<c$ by $(j,\gamma)\map {c\over c-j}\cd\gamma$.  (Read
$(j,\gamma)\map {c!\over c-j}$ if you have taken exponent ${c!\over c-\abs\a}$ in the
definition of coefficient ideals and wish to project to $\N^{n-1}$.)   \>\meds

\<  Let $e_f$ and $e_g$ denote the orders of
the  ideals generated by the equilibrated powers of the coefficients of $f$ and $g$ in $J_P$
and $J_Q$ respectively (not of the coefficient ideals of the ideals generated by $f$ and $g$).
 The equality $e_{f\cd g}/c=\min\{e_f/c_P,e_g/c_Q\}$ then follows from $N(f\cd g)=N(f)+N(g)$
by a computation in $\N^n$.\> \meds

\<Write now elements $h\in K$ as $h= \sum_j a_{h,j}(y)\cd x^j$ so
that $J=\sum_{j<c}(a_{h,j}, h\in K)^{c/c-j}$, and similarly for $P$ and $Q$. There exists
a finite generator set $H$ of $K$ such that  $J=\sum_{j<c}(a_{h,j}, h\in
H)^{c/c-j}$. Let $F$ and $G$ be defined analogously for $P$ and $Q$. Each $h\in H$
is a linear combination of products $f\cd g$ with $f\in F$ and $g\in G$. Enlarging $F$ and
$G$ we may assume that all coefficients in the sum are $0$ or $1$. Replacing $H$ by all
summands of all $h$ we obtain $H=F\cd G$. The coefficients of the elements of the new $H$
generate again $J$.  By the formula above we get \>\meds \<

\hskip 2.4cm $e/c=\min_{h\in H}\,e_h/c=$ \>\meds \<

\hskip 3cm $=\min_{fg\in H}\,\min\,\{e_f/c_P,e_g/c_Q\}=$ \>\meds \<

\hskip 3cm $=\min\,\{\min_{f\in F}\,e_f/c_P,\, \min_{g\in G}\,e_g/c_Q\}=$ \>\meds \<

\hskip 3cm $=\min\,\{e_P/c_P,e_Q/c_Q\}$.
\>\meds

\<This proves the claim. Assume now that $V$ has weak maximal contact with $P$, i.e., that
$e_P$ is maximal. If  $e_P/c_P\leq e_Q/c_Q$ then $e/c$ and hence $e$ is already maximal, by 
the above formula. So assume that $e_P/c_P > e_Q/c_Q$, and that $e/c=e_Q/c_Q$ is not
maximal.  Increasing $e_Q/c_Q$ by a coordinate change requires up to permutations a change
$(x,y)\map (x+g(y),y)$ with slope $\ord\,g= e_Q/c_Q$. \>\meds

\<Let $P^\inter$ and $Q^\inter$ be the resulting ideals. Then, as $P\neq 0, 1$ and $e_P/c_P 
> e_Q/c_Q$,  we get  $e_{P^\inter}/c_{P^\inter}=e_Q/c_Q$, hence
$e/c=\min\{e_{P^\inter}/c_{P^\inter}, e_{Q^\inter}/c_{Q^\inter}\}=e_Q/c_Q$ remains constant,
i.e., was already maximal. This shows that  weak maximal contact with $P\neq 0, 1$ implies
weak maximal contact with $K=P\cd Q$.\>\meds 

\bigskip\goodbreak


{\bf Setup} \medskip

Let $\M=(\J,c,D,E)$ be a singular mobile in $W$ with $\J$ a coherent ideal in a locally
closed regular $n$-dimensional subscheme $V$. Write $J_n$ for the stalk of $\J$ at a point
$a$ of $V$.  A {\it punctual setup} ${}^{<7>}$ of $\M$ at $a$ is a sequence $(J_n\to J_1)$ of
stalks of ideals $J_i$ in a local flag $(W_n\to W_1)$ of $V$ at $a$ satisfying for all
$i\leq n$  \medskip

(1) \ $J_i=M_i\cd I_i$ with $M_i=I_{W_i}(D_i\cap W_i)$ and $I_i$ an ideal in $W_i$ at $a$.
\medskip

(2) \ $M_i$ defines a \nc divisor in $W_i$ at $a$. \medskip

(3) \ $W_{i-1}$ has weak maximal contact at $a$ with the composition ideal
$K_i$ in $W_i$ of $(J_i,c_{i+1},D_i,E_i)$. Here, $c_{i+1}$ is the control of
$J_i$ on $W_i$.  It is given for $i<n$ as the order of $K_{i+1}$ in $W_{i+1}$ at
$a$, and $c_{n+1}=c$.  \medskip

(4) \ $J_{i-1}$ is the junior ideal of $K_i$ in $W_{i-1}$.\medskip
 
Setups depend on and are determined by the choice of the local flag  subject to the above
conditions.  They commute with the operations described in {\it equivariance}. \medskip


\<Let $d$ be maximal with $o_d=0$, and set $d=0$ if all $o_i=\ord_aI_i>0$. From (3) follows
that $K_i=P_i\cd Q_i$ for $i>d$, with $P_i=I_i+M_i^{o_i/(c_{i+1}-o_i)}$ and
$Q_i=I_{W_i}(E_i\cap W_i)$. For $i\leq d$ we have $K_i=1$.\>\meds


A {\it tuned setup} of $(\J,c,D,E)$ on an open subscheme $U$ of $V$ is a sequence of
coherent ideal sheaves $\J_n\to \J_1$ in a decreasing flag $W_n\to W_1$ of closed
subschemes of $U$ such that, for any $i\leq n$, the stalks $J_n\to J_i$ at $a$ 
define the truncation of a punctual setup of $(\J,c,D,E)$ for all points $a$ of
$U\cap \ttop(t_n\to t_{i+1})$. Here, $t_i=(o_i,k_i,\mm_i)$ denotes the tag of
$(J_i,c_{i+1},D_i,E_i)$ at $a$.\medskip


\< Thus $o_i=\ord_aI_i$, $k_i=\ord_aK_i$ in $W_i$, and $m_i$ is the combinatorial tag
of the maximal tight shortcut $N_i$ of $M_i$ of order $\geq c_{i+1}$ at $a$. The 
restriction $(J_{i-1}\to J_1)$ is a punctual setup at $a$ in $W_{i-1}$  of the
restricted mobile $(J_{i-1},c_i,D^{i-1},E^{i-1})$ in $W_{i-1}$ where 
$D^{i-1}=(D_{i-1}\to D_1)$ and $E^{i-1}=(E_{i-1}\to E_1)$.\> \meds 

\<Indices indicate, except for handicaps, the dimension of the corresponding ambient scheme. 
The order of an ideal is taken with respect to this ambient scheme. Observe that, by
definition, $c_{i+1}$ is constant on $W_i$ and only defined there, whereas $D_i$ and $E_i$ are
defined  in $W$. More accurately we should write $W_i=W_i(a)$. Then, locally at a point
$b$, we will be able to choose $W_i$ so that $W_i(a)=W_i(b)$ for all $a$ in $\ttop(K_{i+1})$
near $b$.\>\meds

\bigskip\goodbreak


{\bf Invariant} \medskip

Let $\M=(\J,c,D,E)$ be a mobile in $W$ with $\J$ a non zero ideal in $V$. Assume
that $\M$ admits locally on $V$ tuned setups with induced punctual setups $(J_n\to J_1)$.
Set\medskip

\cl{$i_a(\M) =(t_n\to t_1)\in \N^{4n}$}\medskip

with $t_i=(o_i,k_i,\mm_i)$ the tag of $(J_i,c_{i+1},D_i,E_i)$ at $a$ ${}^{<8>}$. 
Equipping $\N^{4n}$ with the lexicographic order this vector satisfies the following
properties.\medskip \goodbreak 

(a) \ $i_a(\M)$ does not depend on the chosen setup of $\M$ at $a$ and
commutes with the operations described in {\it equivariance}.\medskip

(b) \ The map $a\map i_a(\M)$ is upper semicontinuous on $V$. The induced
stratification of $V$ refines the stratification underlying $D$ and $E$.\medskip

(c) \ The top locus $Z$ of $i_a(\M)$ is regular. Locally, $Z$ lies in the top loci of all
$I_i$, $P_i$, $Q_i$ and $K_i$. It only depends on the restriction of  $i_a(\M)$  to
the support of $\J$. \medskip

(d) \ $Z$ is transversal to all $D_i$ and $E_i$.\medskip

The first part of property (a) will be proven in the section {\it ``Independence and
semicontinuity''}, the second part holds by construction of $i_a(\M)$, and property (b)
is proven in {\it ``Transform of mobile''}. The first two parts of property (c) are
proven in {\it ``Top loci''}. The third part holds by construction of $i_a(\M)$, 
since its first component $o_n=\ord_aI_n$ has top locus inside the support of $J$ and
all remaining components defining $Z$ are orders of ideals taken at points of this
locus.  Property (d) is proven in {\it ``Transversality''}. \medskip 


\< Equivariance of the invariant signifies that the invariant remains the same under
smooth morphisms and field extensions. And, if $W\map W^+$ is a closed embedding of $W$
into a regular scheme $W^+$, and $\M^+$ a mobile in $W^+$ inducing by restriction to $W$
the mobile $\M$, then for $a$ in $W$ the restriction of $i_a(\M^+)$ to dimension $n=\dim\
W$ equals the invariant $i_a(\M)$. We will use this property in {\it ``Resolution of
schemes''} only in case the divisors of the handicaps $D$ and $E$ of $\M$ are empty with
the exception of $E_n\subset W$.\>\meds

\<The first two properties hold because they hold for order of
ideals and because mobiles and their setups commute with these operations. The third
property follows from the fact that the order of the ideal $J^+$ in $W^+$ extending
$J$ is $1$ if $\dim\ W<\dim\ W^+$. If $D^+_{n^+}$ and $E^+_{n^+}$ are empty
and $\dim\ W<\dim\ W^+$, the first composite ideal $K^+$ equals $J^+$, and its
coefficient ideal in a hypersurface $V^+$ of $W^+$ is just the restriction of $J^+$ to
$V^+$. This hypersurface can be chosen to contain $W$. By induction on the dimension, we
reach by iteration a coefficient ideal which equals $J$, and the first components of
$i_a(\M^+)$ are either $1$ (as for $o$ and $k$) or $0$ (as for $m$). Therefore the
invariant of $\M^+$ restricts on $W$ to the invariant of $\M$. We omit the reasoning
for arbitrary $D^+$ and $E^+$. \>\meds


\bigskip\goodbreak


{\bf Comments}\medskip

 Mobiles are the minimal resolution datum needed to define at each stage of the 
resolution process the local invariant. They present global information and only depend
on the initial mobile and the transformation laws: $J$ is the ideal we wish to resolve,
i.e., to transform into a principal monomial ideal supported by exceptional components.
The control $c$ indicates how $J$ transforms under blowup; in addition, it prescribes the
objective: to drop the order of the transforms of $J$ below $c$. The
components $D_i$ of $D$ keep track which part of $J_i$ has already been monomialized.
The components $E_i$ of $E$ collect the exceptional components which may fail to be
transversal to $W_{i-1}$. They both live in $W$ and are independent of the choice of the
local flags, whereas the associated ideals $M_i=I_{W_i}(D_i\cap W_i)$ and 
$Q_i=I_{W_i}(E_i\cap W_i)$ live in $W_i$ and depend on the flag.\medskip

Setups are auxiliary local data depending on the flags and uniquely determined by them. 
They are defined so that the resulting invariant does not depend on the chosen flag.
Their tunedness ensures the semicontinuity of the invariant. The passage from $I_i$ to
$P_i$ is necessary because $\ttop(I_i)$ need not be contained in $\ttop(J_i,c_{i+1})$
when the order of $I_i$ has become small. Multiplication of $P_i$ with $Q_i$ allows to
treat the transversality problem simultaneously with the monomialization of $J_i$.
\medskip

The inclusions  $\ttop(I_i)\supset \ttop(P_i)$, $\ttop(E_i)=\ttop(Q_i)\supset \ttop(K_i)$
and \medskip

\hskip 1cm $\ld\supset\ttop(K_{i+1})\supset
\ttop(J_i,c_{i+1})\supset \ttop(P_i)\supset \ttop(K_i) \supset \ld\supset Z$\medskip

imply that the orders of the weak transforms of $I_i$ and $K_i$ do not increase. It is
shown that they can only remain constant in a very specific situation, in which case the
combinatorial tag $m_i$ drops. This establishes the vertical induction. \medskip

Weak maximal contact may not persist under blowup, and osculating hypersurfaces are a
characteristic zero device to achieve this persistence. For technical reasons it is 
appropriate to take $W_i$ osculating with $P_{i+1}$ instead of $K_{i+1}=P_{i+1}\cd
Q_{i+1}$. The key point of the proof is the commutativity of the descent in dimension
via coefficient ideals and the passage to points of the blowup where the order has
remained constant. This allows descending induction on the dimension.   \medskip

\bigskip\goodbreak


{\bf Transform of mobile}\medskip

Suppose given a mobile with empty combinatorial handicap. By horizontal induction on the
dimension it admits locally tuned setups. Thus the invariant is defined. Its top locus is
closed and regular and gives the center of the first blowup.\medskip


Assume then constructed at a certain stage of the resolution process the mobile $\M=(\J,c,D,E)$ in
$W$ with $\J$ an ideal of control $c$ in $V$ and locally tuned setups with invariant
$i(\M)$. If $\J$ has order $<c$ on $V$ the mobile is resolved. If not, let $W'\map W$ be the
blowup of $W$ in the top locus $Z$ of $i(\M)$ with new exceptional component $Y'$. The
transform $\M'=(\J',c',D',E')$ of $\M$ and a locally defined tuned setup of $\M'$ are
constructed as follows ${}^{<9>}$. \medskip\goodbreak


 By the first component of $i_a(\M)$ we will know that the center $Z$ is contained in $V$
thus inducing a blowup $V'\map V$.  If the controlled transform $\J^!$ of $\J$ in $V'$ has
order $<c$ on $V'$ the resolution of $\J$ is completed. If not, set $\J'=\J^!$ and $c'=c$. Let
$a'$ be a point of $V'$ and denote by $J_n'$ the stalk of $\J'$ at $a'$. Assume constructed
for some $i<n$ the truncated mobile ${}^{i+1}\M'$, and, locally on $V'$, flags $W_n'\to
W_{i+1}'$ with truncated tuned setup $\J_n'\to \J_{i+1}'$ of ${}^{i+1}\M'$. \medskip

 We shall define the truncated mobile ${}^i\M'$ and, locally on $W_{i+1}'$, regular
hypersurfaces $W_i'$, such that the junior ideal $J_i'$ of $K_{i+1}'$ in $W_i'$ yields a
punctual truncated setup $J_n'\to J_i'$ of ${}^i\M'$ induced from a tuned truncated setup
$\J_n'\to \J_i'$. The transform $\M'$ of $\M$ is then defined by descending induction on $i$
and will admit locally tuned setups. \medskip

For $a'\in W'$, let $J_n'\to J_{i+1}'$ denote the truncated punctual setup of 
${}^{i+1}\M'$ induced by $\J_n'\to \J_{i+1}'$. Thus $(t_n'\to t_{i+1}')$ is defined.
By {\it ``Independence and semicontinuity''} it will be independent of  the choice of
the flag.  \medskip


For $n\geq j\geq i+1$, let $T_{j+1}'$ be the locus of points in $W'$ where $(t_n'\to
t_{j+1}')$ equals the value of $(t_n\to t_{j+1})$ along $Z$. We agree that $T_{n+1}'=W_n'$.
Let $O_j'$ be the locus in $T_{j+1}'$ where $o_j'$ is positive and equal the value of
$o_j$ along $Z$. The independence of the truncated invariant implies that these loci do not
depend on the chosen flags. By descending horizontal induction on $i$ we shall assume that the
handicaps of ${}^{i+1}\M'$ satisfy for all $j\geq i+1$ the equalities \medskip 

\hskip 2cm  $D_j'=D_j^{*} + (o_j-c_{j+1})\cd Y'$ \hskip 2.2cm  on $T_{j+1}'$, 
\medskip

\hskip 2cm  $D_j'=\emptyset$ \hskip 5.2cm outside $T_{j+1}'$, \medskip

\hskip 2cm  $E_j'=E_j^\weak $ \hskip 4.9cm on  $O_j'$,\medskip

\hskip 2cm  $E_j'=(Y' + \abs E^\weak) - (E_n'+\ld + E_{j+1}')$ \hskip .4cm
outside $O_j'$.\medskip\goodbreak


Here, $D_j^{*}$ denotes the pull-back of $D_j$ under $W'\map W$. Define $T_{i+1}'$ and
$D_i'$ by setting $j=i$ in the above formulas. Then $D_i'$ does not  depend on the
chosen flags. It is stratified with underlying stratification given by $T_{i+1}'$. 
This establishes the second part of property (b) of the invariant for $D_i'$.  It is
shown in {\it ``Transversality''} that $D_i'$ is a \nc divisor in $W'$. It is labelled
as follows. Shortcuts which do not involve $Y'$ get the label of their image in $W$. The
remaining shortcuts are labelled arbitrarily by distinct and pairwise different
numbers.\medskip


\<The coherence of $D_i'$ along $T_{i+1}'$ is proven by vertical induction. Assume that $D_i$
is coherent along $T_{i+1}$. If $T_{i+1}'=\emptyset$, nothing is to show. If not, $T_{i+1}'=
\ttop(t_n'\to t_{i+1}')$ lies  over $\ttop(t_n\to t_{i+1})$, since the truncated
invariant does not increase under blowup, see {\it ``Decrease of invariant''}.  This and the
coherence of $D_i$ on $\ttop(t_n\to t_{i+1})$ imply that $D_i^*$ is coherent on $T_{i+1}'$.
As both $o_i$ and $c_{i+1}$ are constant on $Z$, the factor $I(Y')^{ o_i-c_{i+1}}$ is coherent
on $Y'$ and hence on $W'$. Therefore $D_i'$ is coherent on $T_{i+1}'$. \>\meds


Let $K_{i+1}'$ be the composition ideal of $(J_{i+1}',c_{i+2}',D_{i+1}',E_{i+1}')$ in
$W_{i+1}'$. It is stratified in $W_{i+1}'$ with strata given by the order $o_{i+1}'$
of $I_{i+1}'$ and the stratifications underlying $D_{i+1}'$ and $E_{i+1}'$. Choose in a
neighborhood of $a'$ a hypersurface $W_i'$ in $W_{i+1}'$ which has weak maximal contact
with $K_{i+1}'$ along $\ttop(K_{i+1}')$.  \medskip

Define $J_i'$ as the junior ideal of $K_{i+1}'$ in $W_i'$. In {\it ``Commutativity''} it is
shown that $W_i'$ can be taken osculating for $P_{i+1}'$, and that then $J_n'\to J_i'$ form a
truncated punctual setup of ${}^i\M'$ along $\ttop(t_n'\to t_{i+1}')$. In particular,
$J_i'=M_i'\cd I_i'$ with $M_i'=I_{W_i'}(D_i'\cap W_i')$ a \nc divisor in $W_i'$, and
$o_i'=\ord_{a'}I_i'$ is defined for $a'$ in $W_i'$.  Let $O_i'$ be the locus of points $a'$ in
$T_{i+1}'$ where $o_i'$ is positive and equal the value of $o_i$ along $Z$. In {\it
``Independence and semicontinuity''} it is shown that $o_i'$ and hence $O_i'$ do not depend on
the choice of $W_i'$. Define \medskip

\hskip 2cm  $E_i'=E_i^\weak$ \hskip 4.9cm on  $O_i'$,\medskip

\hskip 2cm  $E_i'=(Y'+\abs E^\weak)- (E_n' + \ld + E_{i+1}')$ \hskip .5cm
outside $O_i'$.\medskip\goodbreak


\< The divisor $E_i'$ shall collect the exceptional components which may not be
transversal to $W_{i-1}'$. As $W_{i-1}'$ will equal on $O_i'$ the transform of $W_{i-1}$
the new exceptional component $Y'$ will be transversal to $W_{i-1}'$, so that we set
$E_i'=E_i^\weak$ on $O_i'$. Outside, a new $W_{i-1}'$ will be chosen, and hence need not
be transversal to any exceptional component. Thus $E_i'$ consists of all exceptional
components except those taken care of in $E_n' + \ld + E_{i+1}'$. Observe that
$Z\subset\ttop(P_i)\cap\ttop(Q_i)$ implies that $Z$ is contained for $i>d$ in the
components of $E_i$ it meets, and is hence transversal to all components of
$E_i$. \>\meds


It is shown in {\it ``Transversality''} that $E_i'$ is a  reduced  \nc divisor in 
$W'$. It does not depend on the flag $W_n'\to W_i'$.  By vertical induction, we may
assume that $E_i$ is coherent along $O_j\setminus O_{j-1}$ for all $j\geq i$. As $O_j'=
\ttop(t_n'\to t_{j+1}', o_j')$ lies over $\ttop(t_n\to t_{j+1},o_j)$ we conclude that
$E_i'$ is coherent along $O_j'\setminus O_{j-1}'$ for all $j\geq i$, proving for $E_i'$
the second part of property (b) of the invariant.   This completes the construction of
${}^i\M'$ and of its locally tuned truncated setups $\J_n'\to \J_i'$. The transform
$\M'$ of $\M$ and $\J_n'\to \J_1'$ are then defined by descending horizontal induction
on $i$. \medskip

\bigskip\goodbreak


\med \cl {\bf  PROOFS}\bigskip

{\bf Logical structure}\medskip

 In {\it ``Top loci''} it is shown by descending horizontal induction on the dimension that,
locally at points $a$ in $W$, the center $Z$ is contained in the members $W_i$ of the local
flag at $a$ for all $i\geq d$, where $d$ is maximal with $o_d=0$. \medskip

In {\it ``Commutativity''} it is shown, assuming that the assertions of the later section
{\it ``Transversality''} hold at the present stage $W$ of the resolution process, that
at points of $W'$ where a truncation of the invariant at a certain index $i$ has remained
constant, the subsequent descent in dimension and the truncation at the next index
$i-1$ commute with blowup. This and the inclusions of {\it ``Top loci''} allow to show
in {\it ``Decrease of invariant''} that the complete invariant cannot increase when
passing from $W$ to $W'$. From horizontal induction on the dimension then follows that it
actually decreases. This in turn is used together with {\it ``Commutativity''} and {\it
``Top loci''} to show in {\it ``Transversality''} that the handicaps in $W'$ are normal
crossings divisors and transversal to the next center. In this way, the circle of
implications winds up like a spiral through the resolution process.  The sections {\it
``Independence and semicontinuity''} and {\it ``Order of coefficient ideals''} show that
the invariant does not depend on the local choices and is upper semicontinuous. 
\medskip 

\bigskip\goodbreak


{\bf Top loci}\medskip

Let $(J_n\to J_1)$ be a punctual setup of $\M$ at $a$. Let $d$ be maximal with $o_d=0$. We show
that $\ttop(t_n\to t_{d+1})$ coincides with $\ttop(t_{d+1})= \ttop(o_{d+1},
k_{d+1})=\ttop(K_{d+1})$ in $W_{d+1}$. Locally at $a$ this locus lies in $W_d$ if $d\geq 1$,
and equals $a$ if $d=0$. \medskip

Assume that $\ttop(t_n\to t_{i+1})=\ttop(t_{i+1})=\ttop(K_{i+1})$ holds for some $i>d$. Then
$\ttop(t_n\to t_{i+1},t_i)= \ttop(t_{i+1},t_i)= \ttop(t_i)$ in $W_i$ because
$\ttop(t_{i+1})=\ttop(K_{i+1})$ contains $\ttop(J_i,c_{i+1})\supset\ttop(P_i)$ and hence
$\ttop(t_i)$. Descending horizontal induction then yields the assertion.  \medskip  

We show that $Z=\ttop(i(\M))$ equals, locally at $a$, the top locus in $W_d$ of the maximal
tight shortcut $N_d$ of $M_d$ of order $\geq c_{d+1}$ if $d\geq 1$, and $Z=\{a\}$ if
$d=0$. In the second case, $K_1$ is bold regular and different from $1$ in $W_1$ with support
$Z=\{a\}$. If $d\geq 1$ then $o_d=0$ implies $K_d=1$, $k_d=0$, so that $t_d=(0,0,\mm_d)$ and
the remaining components of the invariant are zero. By the above $Z=\ttop(i(\M))$ lies in
$W_d$. It hence  equals $\ttop(\mm_d)=\ttop(N_d)$ in $W_d$.\medskip 

\bigskip\goodbreak


{\bf Commutativity}\medskip

 Suppose given a mobile $\M=(\J,c,D,E)$ in $W$ at a certain stage of the resolution process,
with truncated transform ${}^{i+1}\M'$ in $W'$ as defined in {\it ``Transform of mobile''}, for
some $i<n$.\medskip

Assume constructed truncated punctual setups $(J_n\to J_{i+1})$ of ${}^{i+1}\M$ in local
flags $W_n\to W_{i+1}$ at $a$ and $(J_n'\to J_{i+1}')$ of ${}^{i+1}\M'$ in local
flags $W_n'\to W_{i+1}'$ at $a'$ such that $W_j'=W_j^\weak$ on $O_{j+1}'$, $J_j'=J_j^!$ and
$I_j'=I_j^\weak$ on $T_{j+1}'$, and $P_j'=P_j^\weak$,  $Q_j'=Q_j^\weak$ and thus
$K_j'=K_j^\weak$ on $O_j'$ for all $j\geq i+1$. In particular, $D_i'$ can then be defined as
in {\it ``Transform of mobile''}. Assume also that $W_j$ is chosen osculating for
$P_{j+1}$ at $a$ if $j\geq i$, and that $W_j'$ is chosen osculating for $P_{j+1}'$ at
$a'$ if $j\geq i+1$. \medskip

We show that there exists, locally at $a'$ in $W_{i+1}'$, a regular hypersurface 
$W_i'$ which is osculating for $P_{i+1}'$ and such that the above commutativity
relations also hold for $j=i$, where the ideals $J_i'$, $I_i'$, $P_i'$, $Q_i'$ and
$K_i'$ are defined as in {\it ``Transform of mobile''}. This allows in particular to
define $E_i'$ and the truncated mobile ${}^i\M'$ as in {\it ``Transform of mobile''}.
Moreover, $(J_n'\to J_i')$ will define a truncated punctual setup of ${}^i\M'$ at $a'$,
completing the induction step from $i+1$ to $i$.  \medskip


\<For $W_j'=W_j^\weak$ we need that $P_{j+1}'=P_{j+1}^\weak$ and $o_{j+1}'=
o_{j+1}$ on $O_{j+1}'$. For $J_j'=J_j^!$ we need that $K_{j+1}'=K_{j+1}^\weak$ and
$k_{j+1}'= k_{j+1}$ on $T_{j+1}'$. For $I_j'=I_j^\weak$ we need that $K_{j+1}'=K_{j+1}^\weak$ and
$k_{j+1}'= k_{j+1}$ and $D_j'=D_j^{*} + (o_j-c_{j+1})\cd Y'$ on $T_{j+1}'$. For
$P_j'=P_j^\weak$ we need that $I_j'=I_j^\weak$ and $o_j'=o_j$ on $O_j'$. For $Q_j'=Q_j^\weak$ we
need that $E_j'=E_j^\weak$ on $O_j'$. For $K_j'=K_j^\weak$ we need that $P_j'=P_j^\weak$ and
$Q_j'=Q_j^\weak$ on $O_j'$ for $j>d'=d$, respectively $o_j'=o_j=0$ on $O_j'$ for $j\leq
d'$.\>\meds

\<As  $Z\subset\ttop(P_j)\cap\ttop(Q_j)$ for $j>d$ and hence $o_j=\ord_aP_j=\ord_ZP_j$ and
$r_j=\ord_aQ_j=\ord_ZQ_j$ we get on $O_j'$ from $P_j'=P_j^\weak$ and $Q_j'=Q_j^\weak$ for $j>d'$
that\medskip

\hskip 2cm $K_j'=P_j'\cd Q_j'=P_j^\weak\cd Q_j^\weak =$\medskip

\hskip 2.5cm $=P_j^*\cd Q_j^*\cd I_{W_j'}(Y'\cap W_j')^{-(o_j+r_j)}=$\medskip

\hskip 2.5cm $=(P_j\cd Q_j)^\weak=K_j^\weak$.\medskip 

Observe that we have $I_j'=I_j^\weak$ on $T_{j+1}'$, whereas $P_j'=P_j^\weak$ holds only on
$O_j'$, by the very definition of companion ideals. \> \meds


As $W_i$ is osculating for $P_{i+1}$ it has weak maximal contact with $P_{i+1}$ and hence
with $K_{i+1}$. As $P_{i+1}'=P_{i+1}^\weak$ holds on $O_{i+1}'$ and has order
$o_{i+1}'=o_{i+1}$ there, $W_i^\weak$ is osculating for $P_{i+1}'$, hence has weak maximal
contact with $K_{i+1}'$ on $O_{i+1}'$. We set $W_i'=W_i^\weak$ on $O_{i+1}'$. Outside 
$O_{i+1}'$ we take locally along $\ttop(P_{i+1}')$ for $W_i'$ any hypersurface in $W_{i+1}'$
which is osculating for $P_{i+1}'$. Let $J_i$ and $J_i'$ be the junior ideals of $K_{i+1}$ and
$K_{i+1}'$ in $W_i$ and $W_i'$ respectively.\medskip


\<We use here that points where the order of $P_{i+1}'$ has remained constant lie in
$W_i^\weak$. Note that the local choices of $W_i'$ need not patch globally along
$\ttop(K_{i+1}')$. \>\meds


As the descent to coefficient ideals commutes with weak, respectively controlled transforms,
we get by definition of junior ideals that $J_i'=J_i^!$ on $T_{i+1}'$.\medskip


\< We use here that by {\it ``Top loci''}, the center $Z$ is contained in $W_i$ 
locally at points of $W$ for all $i\geq d$, so that transforms of ideals in $W_i$ are
well defined. \>\meds

\<If both $K_{i+1}$ and $K_{i+1}'$ are not bold regular or $1$, the
equality of coefficient ideals implies that $J_i'=J_i^!$ on  $T_{i+1}'$. If $K_{i+1}'$
is bold regular or $1$, the equality implies the same for $K_{i+1}$, because coefficient
ideals in a hypersurface of weak maximal contact are zero if and only if the ideal is
bold regular or $1$. If $K_{i+1}=1$ then $K_{i+1}'=K_{i+1}^\weak=1$ on $T_{i+1}'$  and
$J_i'=J_i^!=1$. If $K_{i+1}$ is bold regular and different $1$, then $o_{i+1}>0$ by 
definition of $K_{i+1}$. Therefore $\mm_{i+1}=0$ and $i_a(\M_i)=i_a(1)=0$  so that $Z$
equals the support of $K_{i+1}$. Here $\M_i$ denotes the mobile $(J_i, c_{i+1}, D_i\to
D_1, E_i\to E_1)$ in $W$. Then $K_{i+1}'$ equals $1$ and $a'$ lies outside $T_{i+1}'$.
This shows that in all cases $J_i'=J_i^!$ on $T_{i+1}'$.\>\meds


The definition of $M_i'=I_{W_i'}(D_i'\cap W_i')$ with $D_i'$ as in {\it ``Transform of
mobile''} implies that $J_i'=M_i'\cd I_i'$ with $I_i'=I_i^\weak$ on $T_{i+1}'$, and $I_i'=J_i'$
outside. Thus property (1) of setups holds for $J_n'\to J_i'$, and (3) and (4) follow from the 
construction.  As for property (2), we may assume by the assertions of {\it
``Transversality''} applied to $W$ that $W_i$ and $Z$ are transversal to $D_i$. As
$W_i'=W_i^\weak$ on $O_{i+1}'$ and $D_i'=\emptyset$ outside $T_{i+1}'\subset O_{i+1}'$
we conclude that $W_i'$ is transversal to $D_i'$. Hence $M_i'$ is a normal crossings
divisor in $W_i'$ at $a'$. \medskip

The truncated setup $J_n'\to J_i'$ is by construction induced by a locally tuned truncated
setup of ${}^i\M'$ along $\ttop(t_n'\to t_{i+1}')$. By definition of $P_i'$ and $E_i'$ we have 
$P_i'=P_i^\weak$ and $Q_i'=Q_i^\weak$ and thus $K_i'=K_i^\weak$ on $O_i'$. This completes the
induction step from $i+1$ to $i$.\medskip


\<If $J=M\cd I$ in $V$ with controlled transform $J^!=I(Y'\cap V')^{-c}\cd J^*$ in $V'$ and 
if $M'=I(Y'\cap V')^{o-c}\cd M^*$ with $o=\ord_aI$ and $I'=I^\weak=I(Y'\cap V')^{-o}\cd
I^*$, we get\>\meds

\<\hskip 2.2cm  $J'=I(Y'\cap V')^{-c}\cd J^*=$\medskip

\hskip 2.6cm  $=I(Y'\cap V')^{-c}\cd (M\cd I)^*=$\medskip

\hskip 2.6cm $=I(Y'\cap V')^{-c}\cd M^*\cd I^*=$\medskip

\hskip 2.6cm  $=I(Y'\cap V')^{-c}\cd I(Y'\cap V')^{c-o}\cd M'\cd I(Y'\cap
V')^o\cd I=$\medskip

\hskip 2.6cm $=M'\cd I'$.\>\meds


\<We have the following formulas on $O_i'$.  As $M_i'= M_i^*\cd I_{W_i'}(Y'\cap
W_i')^{o_i-c_{i+1}}$ and $c_{i+1}'=c_{i+1}$ and $o_i'=o_i$ we get\meds

\hskip 2cm $P_i^\weak =(I_i+M_i^{o_i/c_{i+1}-o_i})^\weak=$\medskip

\hskip 2.6cm $
=I_i^\weak+(M_i')^{o_i/c_{i+1}-o_i}=$\medskip

\hskip 2.6cm $
=I_i^\weak+(M_i')^{o_i/c_{i+1}-o_i}=$\medskip

\hskip 2.6cm $=I_i'+(M_i')^{o_i'/c_{i+1}'-o_i'}=$\medskip

\hskip 2.6cm $=P_i'$.\meds
  
\>\<On the other hand, $P_{i+1}'=P_{i+1}^\weak$ on $O_{i+1}'$ and thus $W_i'=W_i^\weak$.
As $E_i'=E_i^\weak$ on $O_i'$ we get, setting $r_i=\ord_aQ_i$, \>\meds \<

\hskip 2cm $Q_i^\weak=(I_{W_i}(E_i\cap W_i))^\weak=$\>\meds \<

\hskip 2.6cm $=(I_{W_i}(E_i\cap W_i))^*\cd I_{W_i'}(Y'\cap W_i')^{-r_i}= $\>\meds \<

\hskip 2.6cm  $=(I_{W_i'}(E_i^*\cap W_i')\cd I_{W_i'}(Y'\cap W_i')^{-r_i}=$\>\meds \<

\hskip 2.6cm $=I_{W_i'}(E_i'\cap W_i')=$\>\meds \<

\hskip 2.6cm $=Q_i'$.\>\meds \<  

The superscript $^*$ in the second line above denotes the total transform of
$I_{W_i}(E_i\cap W_i)$ under the blowup $W_i'\map W_i$, whereas in the third line it denotes
the total transform of $E_i$ under $W'\map W$. \medskip

\>\<Observe that if $o_i'=\ord_{a'}I_i'=\ord_{a'}I_i^\weak$ is smaller than $o_i$ and if
$M_i$ appears in $P_i$ (i.e., if  $0<o_i<c_{i+1}$) then $M_i'=M_i^*\cd I_{W_i'}(Y'\cap
W_i')^{o_i-c_{i+1}}\neq M_i^*\cd I_{W_i'}(Y'\cap W_i')^{o_i'-c_{i+1}'}$, though still
$M_i'=I_{W_i'}(D_i'\cap W_i')$ and $J_i'=M_i'\cd I_i'$. Hence $P_i'$ and $K_i'$ are in this
case not the weak transforms of $P_i$ and $K_i$. In particular, $W_{i-1}'$ will not be the
weak transform of $W_{i-1}$. Hence  $E_i'$ must then equal the whole exceptional locus
minus $E_n'+\ld + E_{i+1}'$  to have the center transversal to the exceptional
locus.\>\meds


\bigskip\goodbreak


{\bf Decrease of invariant}\medskip

Let $\M'=(\J',c',D',E')$ be the transform of $\M=(\J,c,D,E)$ by blowing up
$Z=\ttop(i(\M))$ in $W$. We show that $i_{a'}(\M')<i_a(\M)$ for $a\in Z$,
$a'\in Y'$ and $\J\neq 1$ ${}^{<12>}$. \medskip


By {\it ``Top loci''}, the center $Z$ is contained for all $i$ in $\ttop(I_i)$ and
$\ttop(K_i)$. This implies for points $a'$ in $W_i'$ above $a$ that
$\ord_{a'}I_i^\weak\leq \ord_aI_i$ and $\ord_{a'}K_i^\weak\leq \ord_aK_i$. From
{\it ``Commutativity''} follows that $I_i'=I_i^\weak$ and hence $o_i'\leq o_i$ on $T_{i+1}'$, 
and $K_i'=K_i^\weak$ with $k_i'\leq k_i$ on $O_i'$. Combining these inequalities gives
by descending horizontal induction for all points $a'$ of $W'$ and for $d$ maximal with
$o_d=0$ that $(t_n'\to t_{d+1}')\leq (t_n\to t_{d+1})$. \medskip

If $d=0$ or if $d\geq 1$ and $\mm_d=(0,0)$ the center equals locally the support of
$K_{d+1}$. Hence $K_{d+1}'=K_{d+1}^\weak=1$ and $k_{d+1}'=0<k_{d+1}$ on $O_{d+1}'$. If $d\geq
1$ and $\mm_d\neq (0,0)$, a computation in local coordinates using the definition of 
maximal tight shortcuts shows that $\mm_d'<\mm_d$ on $T_{d+1}'$ ${}^{<12>}$. In both cases
we get $i_{a'}(\M')<i_a(\M)$.\medskip

\bigskip\goodbreak 


{\bf Transversality}\medskip

Let $\M=(\J,c,D,E)$ be the mobile obtained at a certain stage of the resolution process for
ideals. We show that $D_i$,  $E_i$  and $\abs E$ are \nc divisors in $W$ transversal to
$Z$, and that $\abs E$ coincides with the exceptional locus accumulated so far. And, if
$W_n\to W_1$ is the local flag of a punctual setup $J_n\to J_1$ of $\M$ at $a$ as constructed
in {\it ``Commutativity''}, $W_i$ is transversal to $D_i$ and $E_i\to E_1$, and not contained 
in $E_i +\ld + E_1$ for all $i$ ${}^{<11>}$. \medskip

 
\<The divisors $E_d\to E_1$ are irrelevant for the definition of $K_d\to K_1$. But it is
necessary to know that $W_d$ is transversal to $E_d\to E_1$ in order to have
$Z=\ttop(N_d)\subset W_d$ transversal to $\abs E$. We have $T_i=\emptyset$ locally at $a$ for
some $i\geq d$ by {\it ``Decrease of invariant''}, hence $O_{d-1}=\emptyset$ if $d\geq 2$ and
then $\abs E= E_n+\ldots + E_{d-1}$.  Note that $E_{d-1}$ may be non empty because $O_d$ need
not be empty when $T_d=\emptyset$.  If $d=1$, then $\abs E= E_n +\ldots + E_1$ by vertical
induction and definition of $E_1$.\>\meds


Assume by vertical induction that these properties hold at the prior stages of the 
process. The definition  of $D_i$ and $E_i$ in {\it ``Transform of mobile''} implies, by
the persistence of normal crossings under blowup in transversal centers, that $D_i$,
$E_i$ and $\abs E$ are \nc divisors. It has been proven in {\it ``Decrease of
invariant''} that some $T_j$ is always empty. If $j\geq 2$ then $O_{j-1}=\emptyset$. The
definition of $E_{j-1}$ if $j\geq 2$ and of $E_1$ if $j=1$ implies by vertical induction
that $\abs E$ coincides with the exceptional locus. We have $\abs E=E_n +\ld + E_j$ and
$E_{j-1}=\ld=E_1=\emptyset$ outside $O_j$. \medskip


We show that $W_i$ is transversal to $D_i$ and $E_i\to E_1$ for all $i$. Let $W\map W^\prior$
denote the last  blowup. By the induction hypothesis, $W_i^\prior$ is transversal to
$D_i^\prior$ and $E_i^\prior \to E_1^\prior$, and $Z^\prior$ is transversal to $D_i^\prior$
and $\abs {E^\prior}$. Recall that $W_i=(W_i^\prior)^\weak$ on $O_{i+1}$, that
$D_i=\emptyset$ outside $T_{i+1}$ and $E_i=\ld =E_1=\emptyset$ outside $O_{i+1}$. This implies
that $W_i$ is transversal to $D_i$ and $E_i\to E_1$.\medskip


\<The transversality of $W_i$ with $D_i$ need not follow from the transversality of $W_i$ with
$E_i$ since $D_i$ could contain other exceptional components than $E_i$, and because $W_i$ will
in general not be transversal to $\abs E$. The transversality of $W_i$ with $D_i$ is
needed in order to know that $M_i$ are \nc divisors in $W_i$.\>\meds 

\<The fact that $W_i$ is transversal to $E_i$ is needed to know that
$\ord_aQ_i$ equals $\ord_aE_i$ and is thus independent of $W_i$. The transversality
of $W_d$ with $E_d$ is needed to know that $Z$ is transversal with $\abs E$,
which in turn implies that $\abs {E'}$ and all $D_i'$ are \nc divisors. Observe that
$W_i$ will in general not be transversal to $E_{i+1}$.\>\meds


We show that $W_i$ is not contained in $E_i +\ld + E_1$. Assume that this holds in
$W^\prior$. Inside $O_{i+1}$ we have $W_i=(W_i^\prior)^\weak$, and $E_i +\ld + E_1$ is
either the weak transform of $E_i^\prior +\ld + E_1^\prior$ or its union with $Y$. Hence,
if $Z^\prior$ is strictly contained in $W_i^\prior$, $W_i$ is not contained in $E_i+\ld+
E_1$ inside $O_{i+1}$. If $Z^\prior=W_i^\prior$ we have $J_i^\prior=1$ because 
$Z^\prior\subset \ttop(J_i^\prior)$ by {\it ``Top loci''}. Therefore $K_{i+1}^\prior$ and hence
$P_{i+1}^\prior$ are bold regular or $1$, with $Z^\prior=W_i^\prior$ the support of
$P_{i+1}^\prior$. Then $o_{i+1}=0<o_{i+1}^\prior$, so we are outside $O_{i+1}$. As $E_i\to
E_1$ are empty there, $W_i$ is not contained in $E_i+\ld+ E_1$ also in this case.\medskip
\goodbreak


\<The fact that $W_i$ is not contained in $E_i$ for $i>d$ is needed to have $Q_i\neq 0$ and
hence $\ttop(K_i)=\ttop(P_i)\cap \ttop(Q_i)$. For $i=d$ it is not needed since
$K_d=1$.\>\meds


We show that $Z$ is transversal to $D_i$, $E_i$ and $\abs E$. As $\abs E$ equals the
exceptional locus the $D_i$ are supported by certain components of $\abs E$. Let $d$ be
the largest index with $o_d=0$ at $a$, setting $d=0$ if all $o_i>0$. By {\it ``Top
loci''}, $Z$ equals locally $\ttop(o_d,k_d,\mm_d)=\ttop(N_d)$ if $d\geq 1$, and $\{a\}$
otherwise.  As $N_d$ is a shortcut of $M_d=I_{W_d}(D_d\cap W_d)$ and $W_d$ is
transversal to $E_d\to E_1$ we see that $Z$ is transversal to $E_d\to E_1$. For $i>d$ we
have $I_i\neq 1$ and $Q_i\neq 0$ so that $\ttop(K_i)=\ttop(P_i)\cap\ttop(Q_i)$. As
$Z\subset \ttop(K_i)$ by {\it ``Top loci''} we get $Z\subset \ttop(P_i)\cap\ttop(E_i)$.
Therefore $Z$ is contained for $i>d$ in all components of $E_i$ it meets. This proves
that $Z$ is transversal to $\abs E$, hence also to all $D_i$ and $E_i$. \medskip

\bigskip\goodbreak


{\bf Independence and semicontinuity}\medskip 

Suppose that we are given, locally on $W$, a truncated tuned setup $\J_n\to \J_i$ of a mobile
$(\J,c,D,E)$ with truncated invariant $(t_n\to t_i)$. We prove by descending horizontal
induction that $(t_n,\ldots,t_i)$ is upper semicontinuous on $W$ and independent of the choice
of the setup. \medskip


\<In case $i=n$, the tag $t_n=(o_n,k_n,\mm_n)$ is constructed without choices in terms 
of the stalk $J_n$ of $\J$, the control $c_{n+1}$ and the divisors $D_n$ and $E_n$. Its
first component $o_n$ is \usc since it is the order of a coherent ideal.  The
definition of $P_n$ and the coherence of $D_n$ imply that the order of $P_n$ is upper
semicontinuous. As also $E_n$ is coherent, the order $k_n$ of $K_n=P_n\cd Q_n$ is upper
semicontinuous. The maximal tight shortcut $N_n$ of $M_n$ has by construction a tag
$\mm_n$ which is upper semicontinuous. Hence $t_n=(o_n,k_n,\mm_n)$ is upper
semicontinuous.  \>\meds


Assume that this holds for $n\to i+1$. Fix a point $b\in W$. Locally, $(t_n,\ldots,t_{i+1})$
attains by semicontinuity at $b$ its maximum. Choose closed
subschemes $W_n,\ldots,W_{i+1}$ of a neighborhood $U$ of $b$ where $\J_n,\ldots,\J_{i+1}$
are defined and induce at all points of $U\cap \ttop(t_n\to t_{i+1})$ a truncated punctual
setup of $\J$.\medskip

After shrinking $U$ there exists a closed hypersurface $W_i$ in $W_{i+1}$ which has weak
maximal contact with $K_{i+1}$ at all points $a$ of $\ttop(K_{i+1})$. By {\it ``Top loci''} we
know that if $o_{i+1}>0$ then $\ttop(t_n,\ldots,t_{i+1})= \ttop(K_{i+1})$. By {\it ``Order
of coefficient ideals'}, $\ord\,J_i$ depends only on the locus $\ttop(K_{i+1})$ and its
transforms under the blowups constructed there. It is hence independent of the choice of the
local flag. The same then holds for $o_i=\ord\,J_i-\ord\,M_i$ because $\ord\,M_i= \ord\,D_i$ by
{\it ``Transversality''}.\medskip

If $o_i>0$, then $K_i=P_i\cd Q_i$ and $k_i=o_i+\ord\,Q_i=o_i+\ord\,E_i$ because $W_i$ is
transversal to $E_i$. If $o_i=0$, then $k_i=0$ and $\mm_i$ only depends on $D_i$ and
$c_{i+1}$. The independence and upper semicontinuity of $t_i=(o_i,k_i,\mm_i)$ follow in both
cases.  \medskip

\bigskip\goodbreak


{\bf Order of coefficient ideals}\medskip

Let $V^\short$ be a closed regular subscheme of a regular scheme $W^\short$ and let $K^\short$ be
a coherent ideal in $W^\short$ with coefficient ideal $J^\short$ in $V^\short$. 
We prove that the order of $J^\short$ in $V^\short$ is determined by the locus
$\ttop(K^\short)\cap  V^\short$ and its transforms under certain blowups ${}^{<13>}$. \medskip

Fix a point $a^\short$ in $\ttop(K^\short)\cap V^\short$ and set
$c=\ord_{a^\short}K^\short$, $e=\ord_{a^\short}J^\short$. Let ${ K}$ and ${ J}$ denote the
ideals generated by $K^\short$ and $J^\short$ in ${ W}=W^\short\times{\Bbb A}^{1}$ and ${
V}=V^\short\times{\Bbb A}^{1}$. Then ${ J}$ is the coefficient ideal of ${ K}$ in ${ V}$. We
have $\ttop({ J},c)\subset\ttop({ K})\cap { V}$. Set ${ L}=\{a^\short\}\times{\Bbb
A}^{1}\subset { V}$ and ${ a}=(a^\short,0)$. Consider the blowup ${ W}'\map{ W}$ of ${ W}$
with center ${ a}$ and exceptional locus ${ Y}'$.  Let ${ K}'$ and ${ J}'$ be the weak and 
controlled transform of ${ K}$ and ${ J}$ with respect to  $c$. From $c'=c$ at all
points of $Y'$ follows that ${ J}'$ is the coefficient ideal of ${ K}'$ in
$V'=V^\weak$, so that $\ttop({ J}',c) \subset \ttop({ K}',c)\cap{ V}'$.\medskip

We have ${ J}'={ I}'\cdot I_{{ V}'}({ Y}'\cap{ V}')^{e-c}$ for some
ideal ${ I}'$ in ${ V}'$. Let ${ L}'$ be the weak transform of ${ L}$. 
Since the order of ${ K}$ along ${ L}$ is $c$, the order of ${ K}'$ along
${ L}'$ is also $c$.  Set ${ a}'={ L}'\cap{ Y}'$.  We blow up ${
W}'$ with center ${ a}'$. Continuing this procedure we get, after $k$ steps, a morphism
${ W}^k\map{ W}^{k-1}$ with exceptional locus ${ Y}^{k}$.  The ideals ${
K}^{k}$ and ${ J}^{k}$ are the weak and controlled transforms of ${
K}^{k-1}$ and ${ J}^{k-1}$ and ${ J}^k$ is the coefficient ideal of ${ K}^k$.
Thus $\ttop({ J}^k,c) \subset\ttop({ K}^k)\cap{ V}^k$. Then the weak transform ${
I}^k$ of $J^{k-1}$  in ${ V}^k$ satisfies ${ J}^k={ I}^k\cdot I_{{ V}^k}({
Y}^k\cap{ V}^k)^{k(e-c)}$. \medskip

Let ${ L}^k$ be the weak transform of ${ L}^{k-1}$. The order of ${ K}^k$
along ${ L}^k$ is $c$. We set ${ a}^k={ L}^k\cap{ Y}^k$. Note that ${
Y}^k\cap{ V}^k\subset\ttop({ J}^k,c)$  if and only if $k(e-c)\geq c$.
For each $k$, let now ${ Y}^k\cap{ V}^k$ be the center of the
next blowup. The associated morphism ${ W}^{k+1}\map{ W}^k$ induces an isomorphism ${
V}^{k+1}\map{ V}^k$ and ${ J}^{k+1}={ J}^k\cdot I({ Y}^k\cap{ V}^k)^{-c}$. These
blowups with center ${ Y}^k\cap{ V}^k$ can be repeated $p_k$ times where $p_k$ is the
integral part of ${k(e-c)\over c}$. Note that $p_k$ depends only on the locus
$\ttop(K^\short)\cap V^\short$ together with its transforms under the previous blowups and on
the control $c$. The assertion then follows from $\ord_a(J^\short)=
(1+\lim\limits_{k\mapsto\infty}{p_k\over k})\cd c$. \medskip 

\bigskip\goodbreak


{\bf Resolution of mobiles} \medskip

Let $\M=(\J,c,D,E)$ be any mobile in $W$ with ideal $\J$ and control $c$ in $V$ which admits
locally on $W$ tuned setups with local flags given by hypersurfaces $W_{i-1}$ in $W_i$
which are osculating for $P_i$. This is e.g. the case when the combinatorial handicap
$D$ is empty, independently of the transversal handicap. We show that $\M$ has a strong
resolution. For this we may assume that $c\leq$ the supremum of the order of $\J$. Blow up
$W$ in $Z=\ttop(i(\M))$ with transformed mobile $\M'=(\J',c',D',E')$ in $W'$. We have seen
in {\it ``Decrease of invariant''} that $i_{a'}(\M')<i_a(\M)$ for $\J\neq 1$ and for all points
$a\in Z$ and $a'\in Y'$. Then $\M'$ admits by {\it ``Commutativity''} locally tuned setups as
before. As the invariant takes values in a well ordered set vertical induction applies.
Hence finitely many blowups make the order of the transform of $\J$ drop below $c$. By {\it
transversality}, all centers are transversal to the exceptional loci. {\it Equivariance}
follows from property (a) of the invariant. This establishes the strong resolution of
mobiles.\medskip

\bigskip\goodbreak


{\bf Resolution of schemes}\medskip

The resolution of mobiles is used to construct a strong resolution of reduced singular
subschemes $X$ of $W$ ${}^{<14>}$. We may assume that $X$ is different from $W$, and
that $W$ is equidimensional. Let $\J$ be the ideal of $X$ in $W$. Associate to it the mobile
$\M=(\J,c,D,E)$ with control $c=1$ and empty handicaps. At any stage $W'$ of the resolution of
$\M$ the controlled transform $\J'$ of $\J$ defines a subscheme of $W'$ formed by the strict
transform $X'$ of $X$ and some components inside the exceptional locus.\medskip

As the final controlled transform of $\J$ equals $1$, there corresponds to each component of $X$
a stage where the strict transform of the component has become regular and has been taken
locally as the center of the next blowup. Let $X_1$ be those components of $X$ which reach this
stage first. The corresponding strict transform $X_1'$ of $X_1$ is regular and, by {\it
``Transversality''} and property (d) of the invariant, transversal to the exceptional
locus. \medskip


Write  $X'=X_1'\cup X_2'$ with $X_2'$ the strict transform of the remaining components 
of $X$. Stop here the resolution process of the mobile $(\J,c,D,E)$ and  define a new
mobile whose resolution implies the separation of $X_2'$ from $X_1'$.  Omitting primes,
let $\K$ be the ideal of $X_2$ in $W$. Let $\J$ be the coefficient ideal of $\K$ in $X_1$
with control $c$ the maximum on $X_1$ of the order of $\K$ in $W$. Set all $D_i$ and
$E_i$ empty with the exception of $E_n$, where $n$ is the dimension of $X_1$, for which
we take the exceptional locus produced so far. \medskip

Resolve the mobile $(\J,c,D,E)$. The controlled transforms of $\J$ are the coefficient 
ideals of the weak transforms of $\K$ as long as the maximum of the order of $\K$ in $W$ along
$X_1$ remains constant. Therefore the resolution of $(\J,c,D,E)$ will make this maximum drop. Hence also
the maximum of the order of the strict transform of $\K$ in $W$ along $X_1$ drops. Iterating
this process the final strict transform $X_2'$ of $X_2$ will be separated from the weak
transform $X_1'$ of $X_1$. Now induction on the number of components applies to construct a
sequence of blowups which makes $X_2'$ regular and transversal to the transversality locus.
Thus $X$ has become regular.\medskip

{\it Embeddedness} and {\it equivariance} follow from the resolution of mobiles. 
The resolution of $X$ does not depend on the embedding of $X$ since, 
under embeddings of $W$ into some $W^+$, the restriction of $i(\M^+)$ to $W$ equals $i(\M)$.
This proves {\it excision}. The construction of a {\it strong resolution} of $X$ is
completed.\medskip 

\bigskip\goodbreak


{\eightrm

{\bf Appendix}\med

We indicate references where concepts and arguments of this paper have their origins or
analogues.  \med 

(1) Maximal tight shortcut: Standard reasoning to prove combinatorial resolution, 
cf. [EV1] 2, p. 114, [BM3] 6.16(b), p. 260. \par

(2) Mobile: Including handicaps extends and refines [H2] 1.1, p. 54, [A1] 3, p. 90, 13, 
p. 208, [EV1] 1.2, p. 112, [BM3] 4.1, p. 241, [BS], p. 407. For the combinatorial
handicap cf. [EV1] 3.1, p. 116, and 4.20, p. 127, [BM3] 4.23, p. 247, for the transversal
handicap cf. [A1] 2, p. 89,  5, p. 98, [EV1] 6.17, p. 149, [BM3] 6.8, p. 256.  \par

(3) Controlled transform: Suitable transform to have coefficient ideals commute with
blowup, cf. [H2] 1.10, p. 57, [EV1] 1.4, p. 113, [BM3] 4.4, p. 242. \par

(4) Composition ideal: Cf. [EV1] 6.20.1, p. 153, 6.21.1, p. 154, 6.21.2, p. 155,
[BM3] 4.23, p. 247. \par

(5) Coefficient ideal: Transfers the resolution problem to smaller dimension, cf.  [H2]
8.5.4, p. 112, [A2] 7.1, p. 19, [EV1] 4.14, p. 122, [BM3] 4.18, p. 246. \par

(6) Weak maximal contact: Characteristic free notion to make the order of the
coefficient ideal coordinate independent, cf. [A2] 1.6, p. 6, [H2], 2.4, p. 63, [AHV] p. 6. 
Osculating hypersurfaces are a characteristic zero device to guarantee weak maximal contact
persistent under blowup, cf. [A1] 13.5, p. 211, [H2] 8.2, p. 106,  8.4, p. 108, [AHV]
1.2.5.7, p. 34, [EV1] 4.4, p. 118,  4.11, p. 121, [BM3] 4.12, p. 244. \par

(7) Setup: References as for composition and coefficient ideal. In contrast to [H2],
1.1, p. 54, 2.3, p. 62, 2.6. and 2.7, p. 67, [EV1] 5.1, p. 129, 5.11, p. 131, [BM3] 4.19, p.
246, patchings and equivalence relations are not needed.  \par

(8) Invariant: Cf. [EV1] 4.16, p. 124, 4.20, p. 127, 6.11, p. 140, 6.13, p. 143 [BM3] 4.20, p.
246, 6.15, p. 259. \par

(9) Transform of mobile: Distinguishes old and new exceptional components, cf. [A1] 2, p. 89,
5, p. 98, [EV1] 1.4, p. 113, 6.17, p. 149, [BM3] 4.4, p. 242, 6.8, p.256.  \par

(10) Commutativity: Needed to construct setups of the transformed mobile, cf. [G] 3.11,
p. 309, [H2] 8.5.6, p. 113, [EV1] 4.6, p. 119, 4.15, p. 123, [BM3] 4.19, p.
246,  4.24, p. 248.\par

(11) Transversality: Cf. [EV1] step 6.2, p. 154, [BM3] 4.12, p. 244.\par

(12) Decrease of invariant: Cf. [A1], 16.5, p. 225, 17.4, p. 233, [H1], p. 312, [EV1] 6.13,
p. 143, [BM3], p. 260.\par

(13) Order of coefficient ideals: The argument is from [H2] 2.8, p. 68, cf. also [A1]
7.4, p. 141.\par

(14) Resolution of schemes: Cf. [EV2] 1.11. \par

}

\med\bigskip\goodbreak

{\bf References}\medskip

{\eightrm \parindent 1.1cm

\litem {[A1]} Abhyankar, S.: Good points of a hypersurface. Adv. Math. { 68} (1988),
87-256.

\litem {[A2]} Abhyankar, S.: Desingularization of plane curves. In: Algebraic Geometry. Arcata
1981, Proc. Symp. Pure Appl. Math. { 40},  Amer. Math. Soc. 

\litem {[AHV]} Aroca, J.-M., Hironaka, H., Vicente, J.-L.: The theory of the maximal
contact. Mem. Mat. Inst. Jorge Juan Madrid { 29} (1975).

\litem {[AJ]} Abramovich, D., de Jong, J.: Smoothness, semi-stability and  toroidal geometry. 
J. Alg. Geometry { 6} (1997), 789-801. 

\litem{[AW]} Abramovich, D., Wang, J.: Equivariant resolution of singularities in characteristic
0.  Math. Res. Letters { 4} (1997), 427-433.

\litem{[BM1]} Bierstone, E., Milman, P.: Uniformization of analytic spaces.  J. Amer. Math. Soc.
{ 2} (1989), 801-836.

\litem{[BM2]} Bierstone, E., Milman, P.: Semianalytic and subanalytic sets. Publ. Math. IHES {
67} (1988), 5-42.
 
\litem{[BM3]} Bierstone, E., Milman, P.: Canonical desingularization in characteristic
zero by blowing up the maxima strata of a local invariant. Invent. math. { 128} (1997),
207-302.

\litem{[BP]} Bogomolov, F., Pantev, T.: Weak Hironaka Theorem. Math. Res. Letters { 3} (1996),
299-307.

\litem{[BS]} Bodn\'{a}r, G., Schicho, J.: Automated Resolution of Singularities for
Hypersurfaces. J. Symb. Comp. 30 (2000), 401-428.

\litem{[BV]} Bravo, A., Villamayor, O.: A strengthening of resolution of
singularities in characteristic zero. To appear in Proc. London Math. Soc. 

\litem{[EV1]} Encinas, S., Villamayor, O.: Good points and constructive resolution of
singularities. Acta Math. { 181} (1998), 109-158.

\litem{[EV2]} Encinas, S., Villamayor, O.: A new theorem of desingularization over fields of
characteristic zero. Preprint.

\litem{[G]} Giraud, J.: Sur la th\'eorie du contact maximal. Math. Z. {
137} (1974), 285-310.

\litem{[Ha1]} Hauser, H.: Resolution of singularities 1860-1999. In: Resolution of
Singularities (ed: Hauser et al.), Progress in Math. { 181}, Birkh\"auser 2000.

\litem{[Ha2]} Hauser, H.: The Hironaka Theorem on Resolution of Singularities.
Forthcoming.

\litem{[H1]} Hironaka, H.: Resolution of singularities of an algebraic variety over a field of
characteristic zero. Ann. Math. { 79} (1964), 109-326.

\litem{[H2]} Hironaka, H.: Idealistic exponents of singularity. In: Algebraic Geometry, The Johns
Hopkins Centennial Lectures 1977, 52-125. 

\litem{[V1]} Villamayor, O.: Constructiveness of Hironaka's resolution. Ann. Sci. \'Ec.
Norm. Sup. Paris { 22} (1989), 1-32.

\litem{[V2]} Villamayor, O.: Patching local uniformizations, Ann. Scient. \'Ec. Norm. Sup.
Paris { 25} (1992), 629-677.

}


{\eightrm

\vskip 1cm Santiago Encinas \hfill Herwig Hauser \par 
Universidad de Valladolid\hfill Universit\"at Innsbruck\par 
E-47014 Spain\hfill A-6020 Austria\par
sencinas@maf.uva.es\hfill herwig.hauser@uibk.ac.at

}


\vfill\eject\end